\documentclass[8 pt]{amsproc}
\oddsidemargin=-1cm\evensidemargin=-1cm
\begin{document}
\numberwithin{equation}{section}

\def\1#1{\overline{#1}}
\def\2#1{\widetilde{#1}}
\def\3#1{\widehat{#1}}
\def\4#1{\mathbb{#1}}
\def\5#1{\frak{#1}}
\def\6#1{{\mathcal{#1}}}

\def\C{{\4C}}
\def\R{{\4R}}
\def\N{{\4N}}
\def\Z{{\4Z}}

\author{Valentin Burcea}
\title{On Huang-Yin's Normal Form}

\begin{abstract}Let  $(z,w)$ be the coordinates in $\mathbb{C}^{2}$. We construct a normal form for a class of real-formal surfaces $M\subset\mathbb{C}^{2}$ defined near a degenerate CR singularity $p=0$ as follows
$$w=z^{2}+\overline{z}^{2} + \mbox{O}\left(3\right).$$
\end{abstract}

\address{V. Burcea: Prejudiced}

\email{valentin@maths.tcd.ie}

\thanks{ THANKS to The Science Foundation Ireland Grant 10/RFP/MTH 2878 for ENORMOUS FUNDING while I was working in Trinity College Dublin, Ireland}

\maketitle



\def\cn{{\C^n}}
\def\cnn{{\C^{n'}}}
\def\ocn{\2{\C^n}}
\def\ocnn{\2{\C^{n'}}}


\def\dist{{\rm dist}}
\def\const{{\rm const}}
\def\rk{{\rm rank\,}}
\def\id{{\sf id}}
\def\tr{{\bf tr\,}}
\def\aut{{\sf aut}}
\def\Aut{{\sf Aut}}
\def\CR{{\rm CR}}
\def\GL{{\sf GL}}
\def\Re{{\sf Re}\,}
\def\Im{{\sf Im}\,}
\def\span{\text{\rm span}}
\def\Diff{{\sf Diff}}

\def\codim{{\rm codim}}
\def\crd{\dim_{{\rm CR}}}
\def\crc{{\rm codim_{CR}}}

\def\phi{\varphi}
\def\eps{\varepsilon}
\def\d{\partial}
\def\a{\alpha}
\def\b{\beta}
\def\g{\gamma}
\def\G{\Gamma}
\def\D{\Delta}
\def\Om{\Omega}
\def\k{\kappa}
\def\l{\lambda}
\def\L{\Lambda}
\def\z{{\bar z}}
\def\w{{\bar w}}
\def\Z{{\1Z}}
\def\t{\tau}
\def\th{\theta}

\emergencystretch15pt \frenchspacing

\newtheorem{Thm}{Theorem}[section]
\newtheorem{Cor}[Thm]{Corollary}
\newtheorem{Pro}[Thm]{Proposition}
\newtheorem{Lem}[Thm]{Lemma}

\theoremstyle{definition}\newtheorem{Def}[Thm]{Definition}

\theoremstyle{remark}
\newtheorem{Rem}[Thm]{Remark}
\newtheorem{Exa}[Thm]{Example}
\newtheorem{Exs}[Thm]{Examples}

\def\bl{\begin{Lem}}
\def\el{\end{Lem}}
\def\bp{\begin{Pro}}
\def\ep{\end{Pro}}
\def\bt{\begin{Thm}}
\def\et{\end{Thm}}
\def\bc{\begin{Cor}}
\def\ec{\end{Cor}}
\def\bd{\begin{Def}}
\def\ed{\end{Def}}
\def\br{\begin{Rem}}
\def\er{\end{Rem}}
\def\be{\begin{Exa}}
\def\ee{\end{Exa}}
\def\bpf{\begin{proof}}
\def\epf{\end{proof}}
\def\ben{\begin{enumerate}}
\def\een{\end{enumerate}}
\def\beq{\begin{equation}}
\def\eeq{\end{equation}}

 \section{Introduction and Main Result}
 
 The study of  Real Submanifolds in Complex Space, near a CR singularity  goes back  to  Bishop\cite{Bi}.  A point $p\in M$  is called a CR
singularity if it is a jumping  discontinuity point for the map $M\ni q\longrightarrow  \dim_{\mathbb{R}}T^{c}_{q}M$  defined near $p$.  Bishop\cite{Bi} considered the case when there exist  coordinates $(z,w)$ in $\mathbb{C}^{2}$ such that near a CR
singularity $p=0$, the surface $M\subset\mathbb{C}^{2}$ is defined locally by
\begin{equation}
 w=z\overline{z}+\lambda\left(z^{2}+\overline{z}^{2}\right)+\rm{O}(3),\label{clasic}
 \end{equation}
where  $\lambda\in\left[0,\infty\right]$ is a holomorphic
invariant called the Bishop invariant. When $\lambda=\infty$, $M$
is understood to be defined by the equation
$w=z^{2}+\overline{z}^{2}+\rm{O}(3)$. If $\lambda$ is
non-exceptional,  Moser-Webster\cite{MW} proved that there
exists a formal transformation that sends $M$ into the following normal form
\begin{equation}
w=z\overline{z}+\left(\lambda+\epsilon
u^{q}\right)\left(z^{2}+\overline{z}^{2}\right),\quad
\epsilon\in\left\{0,-1,+1\right\},\quad q\in\mathbb{N},
\end{equation}
where $w=u+iv$.   Moser\cite{Mos} constructed when $\lambda=0$
the following partial normal form:
\begin{equation}
 w=z\overline{z}+2\Re\left\{\displaystyle\sum_{j\geq s}a_{j}z^{j}\right\}.\label{moser}
\end{equation}

Here $s:=\min\left\{j\in\mathbb{N}^{\star};\hspace{0.1
cm}a_{j}\neq 0\right\}$ is the simplest higher order invariant,
known as the Moser invariant.   When
$s<\infty$, Huang-Yin\cite{HY2}  proved that
\eqref{moser} can be formally transformed into the following normal form 
\begin{equation}
w=z\overline{z}+ 2\Re\left\{\displaystyle\sum_{j\geq
s}a_{j}z^{j}\right\},\quad a_{s}=1,\quad a_{j}=0,\quad
\mbox{if}\quad j=0,1\hspace{0.1 cm}\mbox{mod}\hspace{0.1
cm}s,\quad j>s.\label{hye}
\end{equation}  
 
In this paper, we continue the study of the  C.-R. Singular Real Submanifolds in Complex Spaces considering certain Classes of  Real-Submanifolds   using   motivation from Moser-Webster\cite{MW}. They\cite{MW}   considered   the following class of real-analytic surfaces  
\begin{equation}w=z^{2} +\overline{z}^{2} +\displaystyle\sum _{m+n\geq 3 }a_{m,n}{z}^{m}\overline{z}^{n},\label{00}
\end{equation}
where $(z,w)$ are the coordinates in $\mathbb{C}^{2}$.

Regardless  of its  apparent simplicity, (\ref{00}) defines also a very interesting class of C.-.R. Singular Submanifolds in Complex Spaces.   In particular, it requires a similar approach depending on the Fischer decomposition\cite{Sh} that has been applied by Zaitsev\cite{D1},\cite{D2},\cite{D3} in other situations, and also by the author recently\cite{bu2}.  In order  to develop a partial normal form,  we define the following Fischer-normalization space:

Before beginning, we introduce by \cite{Sh} the following notation
\begin{equation}P^{\star}=\displaystyle\sum_{m+n=k_{0}}\overline{p_{m,n}}\frac{\partial^{m+n}}{\partial z^{m}\partial \overline{z}^{n}},\quad\mbox{if $P(z,\overline{z})=\displaystyle\sum_{m+n=k_{0}}p_{m,n}z^{m}\overline{z}^{n}$.}\label{pol}\end{equation}

In particular, we use the following polynomial
$$Q(z,\overline{z})=z^{2}+\overline{z}^{2}.$$
and in consequence the following differential operator
$$\tr=\frac{\partial^{2} }{\partial z^{2}}+\frac{\partial^{2} }{ \partial\overline{z}^{2}}.$$

 Recalling the Fischer Decomposition from Shapiro\cite{Sh}, we   consider   the following Fischer Decompositions
\begin{equation}   z^{k} =A(z,\overline{z})Q(z,\overline{z})+C(z,\overline{z}),\quad\hspace{0.5 cm}\mbox{where $ \tr \left(C(z,\overline{z})\right)=0$, $\forall k>2$ natural number.} 
\label{V2}
\end{equation} 

We define 
\begin{equation} S_{p},\quad\mbox{for all $p\geq 3$,}
\label{space},\end{equation}
 which consists in real-valued polynomials $P(z,\overline{z})$  of degree $p\geq 1$ in $(z,\overline{z})$ satisfying the normalizations:
$$ P_{k}^{\left(p\right)}(z,\overline{z})=P_{k+1}^{\left(p\right)}(z,\overline{z})Q(z,\overline{z})+R_{k+1}^{\left(p\right)}(z,\overline{z}),\quad \mbox{for all   $k=0,\dots, \left[\frac{p-1}{2}\right]$ and given $P_{0}^{\left(p\right)}(z,\overline{z})=P(z,\overline{z})$,}$$
such that
\begin{equation}
  R_{k+1}^{\left(p\right)}(z,\overline{z})\in   \left(\ker  C^{\star}_{k}  \bigcap   \ker  \overline{C}^{\star}_{k}\bigcap  \ker \tr\right) . 
\end{equation}
 
Furthermore, we assume that 
\begin{equation}W(z,\overline{z})\not\equiv 0,\label{nondeg}
\end{equation}
where we have used the following notation
\begin{equation}\displaystyle\sum _{m+n=3 }a_{m,n}{z}^{m}\overline{z}^{n}\mod \left(C_{3},\overline{C_{3}}\right)W(z,\overline{z}).
\end{equation}

The main result, of this note, is the following

 \bt \label{teo}Let $M\subset\mathbb{C}^{2}$ be a formal surface defined near  $p=0$ by (\ref{00}) satisfying the nondegeneracy condition  (\ref{nondeg}). Then, there exists a unique formal transformation of the following type
\begin{equation}\left(z',w'\right)=\left(z+\displaystyle\sum_{k+l\geq 2}f_{k,l}z^{k}w^{l},\quad w+\displaystyle\sum_{k+l\geq 2}g_{k,l}z^{k}w^{l}\right) ,\label{2800}\end{equation}
that transforms $M$ into the following formal normal form:
\begin{equation}w'=P\left(z',\overline{z'}\right) +\displaystyle\sum _{m+n\geq 3 }a'_{m,n}{z'}^{m}\overline{z'}^{n},\label{3000}
\end{equation}
where  the following Fischer
normalization conditions are satisfied
\begin{equation}\Im \left(\displaystyle\sum _{m+n=p}a'_{m,n}{z'}^{m}\overline{z'}^{n}\right)\in\mathcal{S}_{p-1},\quad \Re\left(\displaystyle\sum _{m+n=p}a'_{m,n}{z'}^{m}\overline{z'}^{n}\right)\in\mathcal{S}_{p},\quad\mbox{for all $p\geq 3$,} \label{cn}
\end{equation}
where $\mathcal{S}_{p}$ is defined in (\ref{space}), and as well the following  normalization conditions holds
\begin{equation} 
W^{\star}\left( R_{3k}^{\left(3\right)}(z,\overline{z})\right)=0,\quad\mbox{for all $k >2$.}
 \label{extra}\end{equation} \et

\textbf{Ackowlodgements} I  acknowledge the importance    of the Grant  06/RFP/MAT018, from Science Foundation of Ireland, in my starting development and  especially in the support  in order to write \cite{V1}, because I did not see this aspect written in the published version of the main part\cite{V1} of my doctoral thesis in Trinity College Dublin. I must emphasize that I did everything which was depending on me in order to improve the writing of the entire components of my doctoral thesis, and  also that the imperfections did not depend on me.  

I apologize to my (former) supervisor Prof. Dmitri Zaitsev for while I was not able to control myself. I have hopes to meet him again, and also for a long and warm friendship, because he is a wonderful person. I empower  any form of funding from Science Foundation of Ireland in Trinity College Dublin. I will not return in Ireland, but I think that   Science Foundation of Ireland must continue to support Mathematics in Trinity College Dublin, with special attention on my (doctoral) supervisor. 
\section{Proof of Theorem \ref{teo}}

\subsection{Notations}Let $(z,w)$ be the holomorphic coordinates in $\mathbb{C}^{2}$. Throughout this note, we   use the following notations
$$a_{\geq l}(z,\overline{z})=\displaystyle\sum_{m+n\geq l}a_{m,n}z^{m}\overline{z}^{n},\quad a_{l}(z,\overline{z})=\displaystyle\sum_{m+n=l}a_{m,n}z^{m}\overline{z}^{n},\quad \mbox{for all $l\geq 3$}.$$

\subsection{  Transformation Equations} Let $M\subset\mathbb{C}^{2}$ be the real-formal surface defined near
$p=0$ by  
\begin{equation}w=Q(z,\overline{z})+\displaystyle\sum_{m+n\geq 3}a_{m,n}z^{m}\overline{z}^{n}.\label{a1}\end{equation}

Let $M'\subset\mathbb{C}^{2}$ be another real-formal   surface defined near $p'=0$ by  
\begin{equation}w'=Q\left(z',\overline{z'}\right)+\displaystyle\sum_{m+n\geq 3}a'_{m,n}{z'}^{m}\overline{z'}^{n}.\label{a2}\end{equation}

We consider $$\left(z',w'\right)= \left(f(z,w),  g(z,w)\right),$$ 
 a formal transformation which sends $M$ into $M'$ and that fixes the point $0\in\mathbb{C}^{2}$. It
follows by (\ref{a2}) that
\begin{equation}g(z,w)=Q\left(f(z,w),\overline{f(z,w)}\right)+\displaystyle\sum_{m+n\geq 3}a'_{m,n}\left(f(z,w)\right)^{m}\overline{\left(f(z,w)\right)^{n}},\label{a3}\end{equation}
where $w$ is defined by (\ref{a1}). Writing as follows
 $$f(z,w)=\displaystyle\sum_{ m+n\geq 0}
f_{m,n}z^{m}w^{n},\quad g(z,w) = \displaystyle\sum_{ m+n\geq 0}
g_{m,n}z^{m}w^{n},$$   it follows by (\ref{a3}) that

\begin{equation}\begin{split}& \displaystyle\sum_{ m+n\geq 0}
g_{m,n}z^{m}\left(Q(z,\overline{z})+a_{\geq 3}(z,\overline{z})\right)^{n}=  Q\left(\displaystyle\sum_{ m+n\geq 0}
f_{m,n}z^{m}\left(Q(z,\overline{z})+a_{\geq 3}(z,\overline{z})\right)^{n},  \displaystyle\sum_{ m+n\geq 0}
f_{m,n}z^{m}\left(Q(z,\overline{z})+a_{\geq 3}(z,\overline{z})\right)^{n}\right)\\&\quad\quad\quad\quad\quad\quad\quad\quad\quad\quad\quad\quad\quad\quad\quad\quad+a'_{\geq 3}\left(\displaystyle\sum_{ m+n\geq 0}
f_{m,n}z^{m}\left(Q(z,\overline{z})+a_{\geq 3}(z,\overline{z})\right)^{n},\displaystyle\sum_{ m+n\geq 0}
f_{m,n}z^{m}\left(Q(z,\overline{z})+a_{\geq 3}(z,\overline{z})\right)^{n}\right). \label{b1}\end{split}\end{equation}

Since our map fixes the point $0\in\mathbb{C}^{2}$, it follows that $g_{0,0}=0$ and $f_{0,0}=0$. Collecting the terms of bidegree $(m, 0)$ in $(z, \overline{z})$ in (\ref{b1}), for all $m<2$, it follows that $g_{m,0}= 0$, for all $m<2$. Collecting the sums of terms of bidegree $(m, n)$ in $(z, \overline{z})$ with $m + n = 2$ in (\ref{b1}), it
follows that
\begin{equation} g_{0,1}Q(z, \overline{z}) = Q \left(f_{1,0}z, \overline{f_{1,0}z}\right).\label{b2} \end{equation}

Then, (\ref{b2}) describes all the possible values of $g_{0,1}$ and $f_{1,0}$ and in particulary 
we obtain that $\Im g_{0,1}=0$. By composing with  n linear
automorphism of the model manifold $\Re w = Q(z, \overline{z})$, we can assume that $g_{0,1}=1$, $f_{1,0}=1$.  By a careful analysis of the terms interactions in (\ref{b1}), we conclude that in order to put suitable normalization conditions, we have to
consider the following terms
$$ g_{m,n}z^{m}\left(P(z, z)\right)^{n} ,\quad f_{m,n}z^{m}Q_{z}(z, \overline{z})\left(Q(z, \overline{z})\right)^{n} ,\quad \overline{f_{m,n}z^{m}Q_{z}(z, \overline{z})}\left(Q(z, \overline{z})\right)^{n}. $$

Collecting the sum of terms of bidegree $(m, n)$ in $(z, \overline{z})$ with $T = m + n$ in (\ref{b1}), it follows that
$$\displaystyle\sum_{m+n=T} \left(a'_{m,n}{z'}^{m}\overline{{z}'}^{n}-a_{m,n}z^{m}\overline{z}^{n}\right)= g_{T}\left (z, Q(z, \overline{z})\right)-2\Re\left\{ Q_{z}(z, \overline{z})f_{T} \left (z, Q(z, \overline{z})\right)\right\} + \dots ,$$
where we have used the following notations
$$ g_{T}(z,w) = \displaystyle\sum_{ m+2n=T}
g_{m,n}z^{m}w^{n},\quad f_{T} (z,w) = \displaystyle\sum_{ m+2n -1=T}
f_{m,n}z^{m}w^{n},$$
and where the terms defined by ''$\dots$,, depend on $f_{k,l}$ with $k + 2l - 1 < T-1$, and as well on $g_{k,l}$ with $k + 2l < T$.

\end{document}